\documentclass{article}
\usepackage{amssymb}
\usepackage{amsmath}

\usepackage{graphicx}

\newtheorem{definition}{Definition}

\newtheorem{lemma}{Lemma}

\newtheorem{theorem}{Theorem}
\numberwithin{equation}{section}

\def\comment#1{{}}

\newcommand{\inte }{{\rm int}\,}

\newcommand{\dom }{\mbox{dom} \,}

\def\qed{{\hfill{\vrule height5pt width3pt depth0pt}\medskip}}

\newcommand{\cover}[1]{\stackrel{#1}{\Longrightarrow}}

\begin{document}

\begin{center}
{\bf \LARGE  Covering Relations  and \\ Non-autonomous
Perturbations of ODEs}

\vskip 0.5cm
 {\large Maciej Capi\'nski}  \\
\vskip 0.2cm { and} \vskip 0.2cm

{\large Piotr Zgliczy\'nski}\footnote{research supported by Polish
State Committee for Scientific
Research grant  2 P03A 041 24
}  \\
\vskip 0.2cm
 Jagiellonian University, Institute of Mathematics, \\
Reymonta 4, 30-059 Krak\'ow, Poland \\ e-mail:
\{mcapinsk,zgliczyn\}@im.uj.edu.pl
\end{center}

\vskip 0.5cm
\begin{center}
\today
\end{center}

\vskip 0.5cm

\begin{abstract}
 Covering relations are a topological tool for detecting periodic
orbits, symbolic dynamics and chaotic behavior for autonomous ODE.
We extend the method of the covering relations onto  systems with
a time dependent perturbation. As an example we apply the method
to non-autonomous perturbations of the R\"{o}ssler equations to
show that for small perturbation they posses symbolic dynamics.
\end{abstract}


\section{Introduction}
\label{sec:intro} The goal of this paper is to answer the
following \newline QUESTION: \emph{Assume that the equation
\begin{equation}
  x'=v(x) \label{eq:auto}
\end{equation}
has a symbolic dynamics (is semiconjugated with some Bernoulli
shift). Consider now small non-autonomous perturbation of
(\ref{eq:auto})
\begin{equation}
  x'(t)=v(x(t)) + \epsilon(t,x(t)). \label{eq:non-auto}
\end{equation}
Will equation (\ref{eq:non-auto}) also have the symbolic dynamics
if $\epsilon$ is sufficiently small?}

We prove that the answer to above question is positive if the
symbolic dynamics is defined in terms of covering relations for
Poincar\'e maps. This is made precise in Section~\ref{sec:appl}
see Theorem~\ref{thm:rossper}. This result is applied to
non-autonomous perturbations  R\"{o}ssler  equations \cite{R} to
show that for small perturbation they posses symbolic dynamics.

The content of the paper can be described as follows. In
Section~\ref{sec:covrel} we recall from paper \cite{GiZ}  the
notion of covering relations for maps. This is the basic technical
tool used in this paper. In Section~\ref{subsec:maintopthm} we
prove Theorem~\ref{thm:basicper} - the basic theorem about
continuation of covering relations for Poincar\'e maps for the
non-autonomous perturbations of ODEs.

In the following sections we apply Theorem~\ref{thm:basicper} to
answer positively our question and present above mentioned
applications.

\section{Topological theorems}

\label{sec:covrel}

\subsection{Covering relations - basic definitions}

\label{subsec:cov-rel}

\begin{definition}
\bigskip \cite{GiZ} An h-set, $N$, is an object consisting of the following
data

\begin{enumerate}
\item  $\left| N\right| $ - a compact subset of $\mathbb{R}^{k}$

\item  $u(N),s(N)\in \{0,1,2,3,\ldots \},$ such that $u(N)+s(N)=k$

\item  a homeomorphism $c_{N}:\mathbb{R}^{k}\rightarrow \mathbb{R}^{k}=%
\mathbb{R}^{u(N)}\times \mathbb{R}^{s(N)}$ such that
\begin{equation*}
c_{N}(\left| N\right| )=\bigskip \overline{B_{u(N)}}(0,1)\times \overline{%
B_{s(N)}}(0,1)
\end{equation*}
\end{enumerate}
\end{definition}

We set
\begin{eqnarray*}
N_{c} &=&\overline{B_{u(N)}}(0,1)\times \overline{B_{s(N)}}(0,1), \\
N_{c}^{-} &=&\partial \overline{B_{u(N)}}(0,1)\times \overline{B_{s(N)}}%
(0,1), \\
N_{c}^{+} &=&\overline{B_{u(N)}}(0,1)\times \partial \overline{B_{s(N)}}%
(0,1), \\
N^{-} &=&c_{N}^{-1}(N_{c}^{-}),\quad N^{+}=c_{N}^{-1}(N_{c}^{+})
\end{eqnarray*}

Later we will quite often drop the parallel lines in $|N|$ and write $N$
instead of $|N|$ to indicate the support of an h-set $N$.

\begin{definition}
\cite{GiZ} \bigskip Assume N,M are h-sets, such that $u(N)=u(M)=u$ and $%
s(N)=s(M)=s$. Let $f:\left| N\right| \rightarrow \mathbb{R}^{k}$
be a continuous map. Let $f_{c}=c_{M}\circ f\circ
c_{N}^{-1}:N_{c}\rightarrow \mathbb{R}^{u}\times \mathbb{R}^{s}$.
We say that
\begin{equation*}
N \cover{f} M
\end{equation*}
($N$ $f$-covers $M$) if the following conditions are satisfied

\begin{enumerate}
\item[1.]  There exists a continuous homotopy $h:[0,1]\times
N_{c}\rightarrow \mathbb{R}^{u}\times \mathbb{R}^{s}$ such that
the following conditions hold true
\begin{eqnarray*}
h_{0} &=&f_{c}, \\
h([0,1],N_{c}^{-})\cap M_{c} &=&\emptyset , \\
h([0,1],N_{c})\cap M_{c}^{+} &=&\emptyset .
\end{eqnarray*}

\item[2.1.]  If $u>0,$ then there exists a linear map $A:\mathbb{R}%
^{u}\rightarrow \mathbb{R}^{u},$ such that
\begin{eqnarray}
h_{1}(p,q) &=&(Ap,0),\text{ \ \ where }p\in \mathbb{R}^{u}\text{
and }q\in
\mathbb{R}^{s},  \label{eq:h1(p,q) = (Ap,0)} \\
A(\partial B_{u}(0,1)) &\subset &\mathbb{R}^{u}\backslash \overline{B_{u}}%
(0,1).  \label{eq:A(brzeg) zawarte w R^u minus B(0,1)}
\end{eqnarray}

\item[2.2.]  If u=0, then
\begin{equation*}
h_{1}(x)=0,\quad \text{for \ }x\in N_{c}.
\end{equation*}
\end{enumerate}
\end{definition}

With above definition we have the following theorem (see also \cite
{MM,Ztmna,ZNon} for its precursors).

\begin{theorem}
\cite{GiZ} \label{thm:t1} Let $N_{i},\,i=0,\ldots ,n$ be an h-set and $%
N_{n}=N_{0}$. Assume that for each $i=1,\ldots ,n$ we have
\begin{equation}
N_{i-1} \cover{f_{i}} N_{i}  \label{eq:nakrywanie}
\end{equation}
then there exists a point $x\in $int$\left| N_{0}\right| $, such
that
\begin{eqnarray*}
f_{i}\circ f_{i-1}\circ \ldots \circ f_{1}(x) &\in
&\text{int}\left|
N_{i}\right| ,\quad i=1,\ldots ,n \\
f_{i}\circ f_{i-1}\circ \ldots \circ f_{1}(x) &=&x.
\end{eqnarray*}
\end{theorem}

\subsection{Continuation of covering relations for Poincar\'e maps for
non-autonomous perturbations}

\label{subsec:maintopthm}

Let $v:\mathbb{R}^{k}\rightarrow \mathbb{R}^{k}$ be a $C^{1}$ function. Let
us consider an autonomous differential equation
\begin{equation}
x^{\prime }=v(x)  \label{eq:x'=v(x)}
\end{equation}

Let $V_{0},V_{2},\ldots ,V_{n-1},V_{n}$ be Poincar\'e sections of the system
generated by the equation (we do not require that they are different). Let $%
1\leq i\leq n$, and let $x$ be the solution of the problem
\begin{eqnarray*}
x^{\prime } &=&v(x) \\
x(0) &=&x_{0}
\end{eqnarray*}
where $x_{0}\in V_{i-1}.$ For $i=1,\ldots ,n,$ by $\sigma_{i}(x_{0})$ we
will denote the first time for which the solution $x$ reaches the section $%
V_{i}$,

\begin{equation*}
\sigma _{i}(x_{0}):=\inf \{t>0:x(t)\in V_{i}\}
\end{equation*}
When it will be evident from the context which sections we wish to consider,
we will sometimes omit the index $i$.

We will also define functions (Poincar\'e maps)
\begin{eqnarray*}
f_{ji} &:&V_{j} \supset \dom(f_i)\rightarrow V_{i} \\
f_{ji}(x_{0}) &:&=x(\sigma _{i}(x_{0})).
\end{eqnarray*}

Now in the context of the question asked in the introduction we assume that
we have a finite set of covering relations for sets $N_i$ for Poincar\'e
maps defined by an ODE, which leads via Theorem~\ref{thm:t1} to symbolic
dynamics. For example assume that we have the following covering relations $%
N_i \cover{P} N_j $ for $i,j=0,1$ (compare topological horseshoes in \cite
{Zbif}). Then from Theorem~\ref{thm:t1} it follows that we have a
semiconjugacy onto the Bernoulli shift on two symbols.

Let us now consider the equation (\ref{eq:x'=v(x)}) with a time dependent
perturbation
\begin{equation}
x^{\prime }(t)=v(x(t))+ \epsilon(t,x(t))  \label{eq:rownanie z
perturbacja}
\end{equation}
We will try to show a similar result (for example topological horseshoe) for
this perturbed equation. It seams very likely that for small perturbation
the above result should hold. Let us start with the fact that for small
perturbations of the equation the covering relations (\ref{eq:nakrywanie})
for the solution still hold. Let us clarify what we will exactly understand
by the functions $f_{ji}$ in the setting of the perturbed equation (\ref
{eq:rownanie z perturbacja}). Let us consider the equation (\ref{eq:rownanie
z perturbacja}) with the following initial conditions
\begin{eqnarray}
x^{\prime } &=&v(x)+\epsilon(t,x)  \notag \\
x(T) &=&x_{0}  \label{eq:problem cauchyego dla rownania z perturbacja} \\
x_{0} &\in &N_{i-1}  \notag
\end{eqnarray}
Let $x$ be the solution of problem (\ref{eq:problem cauchyego dla rownania z
perturbacja}). We will define functions $f_{ji}^{T}$ which will be analogous
to the functions $f_{ji}$. As before
\begin{eqnarray}
f_{ji}^{T} &:&V_{j}\rightarrow V_{i}  \notag \\
f_{ji}^{T}(x_{0}) &:&=x(\sigma _{i}(x_{0},T))
\label{eq:definicja
funkcji fiT} \\
\sigma _{i}(x_{0},T) &:&=\inf \{t>T:x(t)\in V_{i}\}  \notag
\end{eqnarray}
Let us note that for $|\epsilon|$ sufficiently small the functions above are
well defined \cite{GiZ}. What is more if there exists a covering relation $%
N_j \cover{f_{ji}} N_i$ then for small $|\epsilon|$ the term
\begin{equation}
\sigma _{i}(x_{0},T)-T\quad \text{is bounded for all }x_{0}\in N_{j}
\label{eq:sigma-T ograniczone}
\end{equation}
in the sense that if we change the $\epsilon$, then the lower and the upper
bound of the expression does not change for all $x_{0}\in N_{j}$.

The goal of this section is to establish the following

\begin{theorem}
\label{thm:basicper} Let $v:\mathbb{R}^{k}\rightarrow \mathbb{R}^{k}$ be $%
C^1 $-function, let $V_{1},\ldots ,V_{n}$ be the Poincar\'e
sections for the solution of the equation
\begin{equation}
x^{\prime }=v(x) \label{eq:ode-bp}
\end{equation}
Let $N_{i}\subset V_{i}$ , $i=1,\ldots ,n$ be h-sets, we denote
this family by $\mathcal{H}$

Assume that we have a set $\Gamma$ of covering relations $N_i
\cover{f_{ij}} N_j$ for some $N_i, N_j \in \mathcal{H}$, where
$f_{ji}$ are Poincar\'e maps for (\ref{eq:ode-bp}).

Then there exits $\delta = \delta(\Gamma)$ such that for all
continuous $\epsilon :\mathbb{R}^{k+1}\rightarrow R^{k}$ such that
$\left| \epsilon \right| <\delta$ we have:

For any $t_{0}\in \mathbb{R}$ and for any infinite chain of
covering relations from $\Gamma $
\begin{equation*}
N_{0}\cover{f_{01}}N_{1}\cover{f_{12}}N_{2}\cover...
\end{equation*}
where $N_{i}\in \mathcal{H}$ and $\left( N_{i}\cover{f_{i,i+1}}%
N_{i+1}\right) \in \Gamma $ for $i=0,1,\dots $, \newline there
exists a point $x_{0}\in N_{0}$ and a sequence $\left\{
t_{i}\right\} _{i=0}^{\infty }$ , $t_{0}<t_{1}<\ldots
<t_{m}<\ldots $ , such that for the solution $x$ of the equation
\begin{eqnarray}
x^{\prime } &=&v(x)+\epsilon (t,x)  \label{eq:x'=v+e} \\
x(t_{0}) &=&x_{0}  \notag
\end{eqnarray}
we have
\begin{equation*}
x(t_{i})\in \text{int} N_{i},\qquad i=1,2,\dots
\end{equation*}
\end{theorem}

Before we move on to the proof of this theorem we shall need some
preliminary results.

\begin{lemma}
\label{lem:Lemma o nakrywaniu dla perturbacji} Assume that $f_i$
is a Poincar\'e map for (\ref{eq:ode-bp}). If
\begin{equation*}
N_{i-1}\cover{f_{i}}N_{i}
\end{equation*}
then there exists a $\delta >0$ such that for all $\epsilon $ such that $%
\left| \epsilon \right| <\delta $, for all $T\in \mathbb{R}$%
\begin{equation*}
N_{i-1}\cover{f_{i}^{T}}N_{i}.
\end{equation*}
Furthermore for all $i$ there exists a homotopy $H^{i}:[0,1]\times \mathbb{R}%
\times N_{i-1,c}\rightarrow \mathbb{R}^{u}\times \mathbb{R}^{s},$%
\begin{eqnarray}
H^{i}(0,T,x) &=&f_{i,c}^{T}(x),  \label{eq:H1} \\
H^{i}(1,T,(p,q)) &=&(A_{i\,}p,0),  \label{eq:H2} \\
H^{i}([0,1],T,N_{i-1,c}^{-})\cap N_{i,c} &=&\emptyset ,  \label{eq:H3} \\
H^{i}([0,1],T,N_{i-1,c})\cap N_{i,c}^{+} &=&\emptyset .
\label{eq:H4}
\end{eqnarray}
where $x=(p,q)$ and $A_{i\,}$ is the linear map from the
definition of the covering relation for the covering
$N_{i-1}\cover{f_{i}}N_{i}$.
\end{lemma}

This Lemma states that for all the Poincar\'e maps $f_{i}^{T}$ , for any $T$%
, there exists a homotopy $H_{T}^{i}=H^{i}(\cdot ,T,\cdot )$ which
transports the function $f_{i,c}^{T}$ into the linear function $(A_{i\,},0),$
which is independent from $T$. What is more, the family of functions $%
H_{T}^{i}$ is continuous with respect to $T$.

\noindent \textbf{Proof:} The first part of the lemma regarding the fact
that
\begin{equation*}
N_{i-1} \cover{f_{i}^{T}} N_{i}
\end{equation*}
is a consequence of the Theorem 13 from \cite{GiZ}. To prove the second part
of the lemma, let us consider the following differential equation
\begin{equation*}
x^{\prime }=v(x)+(\frac{1}{2}-\lambda )\epsilon(t+T,x)
\end{equation*}
where $\lambda \in \lbrack 0,\frac{1}{2}].$ We can define Poincar\'e maps $%
f_{i}^{\lambda ,T}$ in the same manner as we have defined the functions $%
f_{i}^{T}$ in (\ref{eq:definicja funkcji fiT}). From the first part of the
lemma we know that
\begin{equation*}
N_{i-1}\cover{f_{i}^{\lambda ,T}} N_{i}.
\end{equation*}
Let us note that $f_{i}^{\frac{1}{2},T}=f_{i}$ . Since
\begin{equation*}
N_{i-1}\cover{f_{i}} N_{i}
\end{equation*}
we know that there exists a homotopy $h^{i}:[0,1]\times N_{i-1,c}\rightarrow
\mathbb{R}^{u}\times \mathbb{R}^{s},$ which satisfies the conditions 1, 2.1
and 2.2, from the definition of the covering relation.

We can now define our homotopy as
\begin{equation*}
H^{i}(\lambda ,T,x):=\left\{
\begin{array}{ll}
c_{N_{i}}\circ f_{i}^{\lambda ,T}(x)\circ c_{N_{i-1}}^{-1} & \text{for }%
\lambda \in \lbrack 0,\frac{1}{2}] \\
h^{i}(2\lambda -1,x) & \text{for }\lambda \in (\frac{1}{2},1]
\end{array}
\right.
\end{equation*}

We need to show that this homotopy satisfies the conditions (\ref{eq:H1}), (%
\ref{eq:H2}), (\ref{eq:H3}), (\ref{eq:H4}). The first two conditions are
evident from the definition of $H.$ From the fact that
\begin{equation*}
N_{i-1}\cover{f_{i}^{\lambda ,T}} N_{i}
\end{equation*}
we know that
\begin{eqnarray*}
H^{i}([0,\frac{1}{2}],T,N_{i-1,c}^{-})\cap N_{i,c} &=&\emptyset , \\
H^{i}([0,\frac{1}{2}],T,N_{i-1,c})\cap N_{i,c}^{+} &=&\emptyset .
\end{eqnarray*}
The fact that
\begin{eqnarray*}
H^{i}((\frac{1}{2},1],T,N_{i-1,c}^{-})\cap N_{i,c} &=&\emptyset , \\
H^{i}((\frac{1}{2},1],T,N_{i-1,c})\cap N_{i,c}^{+} &=&\emptyset ,
\end{eqnarray*}
follows from the conditions 2.1 and 2.2 for the covering
\begin{equation*}
N_{i-1}\cover{f_{i}} N_{i}.
\end{equation*}
Hence all the conditions (\ref{eq:H1}), (\ref{eq:H2}), (\ref{eq:H3}), (\ref
{eq:H4}) hold. \qed

The following lemma will be the main tool for the proof of the
Theorem \ref {thm:basicper}.

\begin{lemma}
\label{lem: lemat glowny}Let $v:\mathbb{R}^{k}\rightarrow
\mathbb{R}^{k}$ be $C^{1}$-function, let $V_{1},\ldots ,V_{n}$ be
the Poincar\'{e} sections for the solution of the equation
\begin{equation}
x^{\prime }=v(x) \label{eq:ode-bp2}
\end{equation}
Let $N_{i}\subset V_{i}$ , $i=1,\ldots ,n$ be h-sets. Let $f_i$ be
Poincar\'e maps for (\ref{eq:ode-bp2}).  If
\begin{equation}
N_{0}\cover{f_1}N_{1}\cover{f_2}\dots \cover{f_{n}}N_{n}
\label{eq:chain-cov-rel}
\end{equation}
then there exists a $\delta >0$, $\delta $ depends only on the set
of covering relations in the chain (\ref{eq:chain-cov-rel}) and
not on the
length of the chain, such that for all continuous $\epsilon :\mathbb{R}%
^{k+1}\rightarrow R^{k}$ such that $\left| \epsilon \right|
<\delta $ for all $T\in \mathbb{R}$
\begin{equation*}
N_{i-1}\cover{f_{i}^{T}}N_{i}\text{\quad for }i=1,\ldots ,n
\end{equation*}
and for any $t_{0}\in \mathbb{R}$ there exists a point $x_{0}\in
N_{0}$ and a sequence $t_{0}<t_{1}<\ldots <t_{n}$ , such that for
the solution $x$ of the equation
\begin{eqnarray}
x^{\prime } &=&v(x)+\epsilon (t,x)  \label{eq:x'=v(x) + epsilon(t,x)} \\
x(t_{0}) &=&x_{0}  \notag
\end{eqnarray}
we have
\begin{equation*}
x(t_{i})\in \text{int}N_{i}\quad \text{for }i=1,\ldots ,n
\end{equation*}
\end{lemma}

\textbf{Proof:} From Lemma \ref{lem:Lemma o nakrywaniu dla perturbacji} we
know that the first part of the lemma is true.

\bigskip \noindent Without any loss of generality we will will give the
proof for $t_{0}=0$. We will also assume that
\begin{eqnarray*}
c_{N_{i}} &=&\text{Id }\quad \text{for }i=0,\ldots ,n \\
f_{i} &=&f_{i,c}\quad \text{for }i=1,\ldots ,n \\
\left| N_{i}\right|  &=&N_{i,c}\quad N_{i}^{\pm }=N_{i,c}^{\pm }.
\end{eqnarray*}
Let us define a function
\begin{eqnarray*}
g &:&N_{n}\rightarrow V_{0} \\
g:= &&(A_{n+1},0)
\end{eqnarray*}
where $A_{n+1}:R^{u}\rightarrow R^{u}$ is any linear map such that $%
A_{n+1}(\partial B_{u}(0,1))\subset R^{u}\backslash \overline{B_{u}}(0,1).$
Clearly we have
\begin{equation*}
N_{n}\cover{g}N_{0}
\end{equation*}
This artificial function will be needed to close the loop of
covering relations (compare Thm.~\ref{thm:t1}), so that it is
possible to define the function (\ref{eq: definicja F}) later on.

Let us\ define functions
\begin{eqnarray*}
F_{i} &:&N_{i-1}\times \mathbb{R}\rightarrow V_{i}\times \mathbb{R\quad }%
\text{for }i=1,\ldots ,n \\
F_{i}(x,T):= &&(f_{i}^{T}(x),\sigma _{i}(x,T)).
\end{eqnarray*}
If we start from the set $N_{0},$ then from (\ref{eq:sigma-T ograniczone})
we know that there exists $s_{1},s_{2},\ldots s_{n}$ and $r_{1},r_{2},\ldots
,r_{n}$ such that
\begin{eqnarray}
F_{1}(N_{0},0) &\subset &V_{1}\times \text{int}I_{1}  \notag \\
F_{2}(N_{1},I_{1}) &\subset &V_{2}\times \text{int}I_{2}  \notag \\
&&\ldots   \label{eq:konstrukcja XN} \\
F_{n}(N_{n-1},I_{n-1}) &\subset &V_{n}\times \text{int}I_{n}  \notag
\end{eqnarray}
where
\begin{equation}
I_{j}=[s_{j}-r_{j},s_{j}+r_{j}].  \label{eq:definicja sj+rj}
\end{equation}
Let us define
\begin{equation*}
XN:=(N_{0}\times \lbrack -1,1])\times (N_{1}\times I_{1})\times \ldots
\times (N_{n}\times I_{n})
\end{equation*}
and
\begin{equation}
F:XN\rightarrow (\mathbb{R}^{k}\times \mathbb{R})^{n+1}
\label{eq: definicja F}
\end{equation}
\begin{eqnarray*}
F((x_{0},t_{0}),\ldots ,(x_{n-1},t_{n-1}))=((x_{0}-g(x_{n}) &,&t_{0}), \\
(x_{1}-f_{1}^{t_{0}}(x_{0}) &,&t_{1}-\sigma _{1}(x_{0},t_{0})), \\
&&\ldots , \\
(x_{n}-f_{n}^{t_{n-1}}(x_{n-1}) &,&t_{n}-\sigma _{n}(x_{n-1,}t_{n-1})))
\end{eqnarray*}

We will show that there exists an $x=((x_{0},t_{0}),\ldots
,(x_{n},t_{n}))\in \inte XN$ such that $F(x)=0.$ Once we find the $x,$ we
will have our $x_{0},t_{0},\ldots ,t_{n},$ because from the definition we
know that for $i=1,\ldots ,n$
\begin{equation*}
f_{i}^{T}(x_{0})=x(\sigma _{i}(x_{0},T))
\end{equation*}
and from the fact that $F(x)=0$ we shall have
\begin{eqnarray*}
x_{i}-x(t_{i}) &=&x_{i}-x(\sigma _{i}(x_{i-1},t_{i-1}))\\
 &=& x_{i}-f_{i}^{t_{i-1}}(x_{i-1}) \\
 &=&0
\end{eqnarray*}
which will mean that
\begin{equation*}
x(t_{i})\in \inte N_{i}\quad \text{for }i=1,\ldots ,n
\end{equation*}
What is more, from our construction and the fact that $F(x)=0$ we will know
that $t_{0}=0$ and that
\begin{equation}
t_{i}-\sigma (x_{i-1},t_{i-1})=0  \label{eq:ti-sigma=0}
\end{equation}
which means that
\begin{equation*}
0=t_{0}<t_{1}<\ldots <t_{n}.
\end{equation*}

Our goal is therefore to find the $x\in \inte XN$ for which $F(x)=0.$ Let us
define a homotopy
\begin{equation*}
H:[0,1]\times XN\rightarrow (\mathbb{R}^{k}\times \mathbb{R})^{n+1}
\end{equation*}
\begin{eqnarray*}
H(\lambda ,(x_{0},t_{0}),\ldots ,(x_{n-1},t_{n-1})) &=&((x_{0}-G(\lambda
,t_{n},x_{n}),t_{0}), \\
&&(x_{1}-H^{1}(\lambda ,t_{0},x_{0}),t_{1}-\lambda s_{1}-(1-\lambda )\sigma
(x_{0},t_{0})), \\
&&, \ldots , \\
&&(x_{n}-H^{n}(\lambda ,t_{n-1},x_{n-1}),t_{n}-\lambda s_{n} \\
&&-(1-\lambda )\sigma (x_{n-1},t_{n-1})))
\end{eqnarray*}
where for $i=1,\ldots ,n$, $H^{i}$ is the homotopy from the Lemma \ref
{lem:Lemma o nakrywaniu dla perturbacji} and $G(\lambda ,\cdot ,\cdot
)=(A_{n+1},0)$ . Let us note that $H(0,x)=F(x)$ and that
\begin{equation}
H(1,x)=B(x-((0,0),(0,s_{1}),\ldots ,(0,s_{n})))  \label{eq:H(1,.)=B}
\end{equation}
where
\begin{eqnarray*}
B((x_{0},t_{0}),\ldots ,(x_{n-1},t_{n-1}))
&=&(((p_{0},q_{0})-(A_{n+1}p_{n},0),t_{0}), \\
&&((p_{1},q_{1})-(A_{1\,}p_{0},0),t_{1}), \\
&&\ldots , \\
&&((p_{n},q_{n})-(A_{n\,}p_{n-1},0),t_{n}))
\end{eqnarray*}
\begin{equation*}
x_{i}=(p_{i},q_{i})\quad \text{for }i=0,\ldots ,n
\end{equation*}

Let us assume that we have the following two lemmas which we will prove
after completing this proof

\begin{lemma}
\label{lem:lemat dla lambda=1}
\begin{equation*}
\text{deg}(H(1,\cdot ),\text{int}XN,0)=\pm 1
\end{equation*}
\end{lemma}

\begin{lemma}
\label{lem:lemat dla lambda in [0,1]}For all $\lambda \in \lbrack
0,1]$ the local Brouwer degree deg$(H(\lambda ,\cdot ),$int$XN,0)$
is defined, constant and independent from $\lambda .$
\end{lemma}

Let us now complete our proof using the two lemmas. From Lemmas \ref
{lem:lemat dla lambda=1} and \ref{lem:lemat dla lambda in [0,1]} we know
that
\begin{equation*}
\deg (F,\text{int}XN,0)=\deg (H(0,\cdot ),\text{int}XN,0)=\deg (H(1,\cdot ),%
\text{int}XN,0)=\pm 1
\end{equation*}
which means that there exists an $x \in \text{int}XN $ such that
$F(x)=0$,
\newline
$x=((x_{0},t_{0}),\ldots ,(x_{n},t_{n}))$ hence we have found our $x_{i}\in %
\inte N_{i}$ and $t_{i}$. \qed

\bigskip

Now to finish of the argument, let us prove the Lemmas \ref{lem:lemat dla
lambda=1} and \ref{lem:lemat dla lambda in [0,1]}.

\noindent \textbf{Proof of Lemma \ref{lem:lemat dla lambda=1}:} From (\ref
{eq:H(1,.)=B}) we know that
\begin{equation*}
H(1,x)=B(x-((0,0),(0,s_{1}),\ldots ,(0,s_{n})))
\end{equation*}
where $B$ is linear. From the degree for affine maps (\ref{eq:deg(f,C)
=sgn(det(B))}) we have
\begin{equation*}
\deg (H(1,\cdot ),\text{int}XN,0)=\text{sgn}(\text{det}B)
\end{equation*}

which means that to prove the lemma it is sufficient to show that $B$ is an
isomorphism. Let us recall the definition of $B.$%
\begin{eqnarray*}
B((x_{0},t_{0}),\ldots ,(x_{n-1},t_{n-1}))
&=&(((p_{0},q_{0})-(A_{n+1}p_{n},0),t_{0}), \\
&&((p_{1},q_{1})-(A_{1\,}p_{0},0),t_{1}), \\
&&\ldots , \\
&&((p_{n},q_{n})-(A_{n\,}p_{n-1},0),t_{n}))
\end{eqnarray*}
\begin{equation*}
x_{i}=(p_{i},q_{i})\quad \text{for }i=0,\ldots ,n-1
\end{equation*}

We have to show that $B(x)=0$ implies $x=0$. If $B(x)=0$ then
\begin{eqnarray*}
t_{0} &=&t_{1}=\ldots =t_{n}=0 \\
q_{0} &=&q_{1}=\ldots =q_{n}=0.
\end{eqnarray*}
We also know that
\begin{eqnarray}
p_{0} &=&A_{n+1}p_{n}  \notag \\
p_{1} &=&A_{1}p_{0}  \notag \\
&&\ldots  \label{eq:pi=Api-1} \\
p_{n} &=&A_{n}p_{n-1}  \notag
\end{eqnarray}
which means that
\begin{equation*}
p_{0}=A_{n+1}\circ \ldots \circ A_{1}p_{0}
\end{equation*}

The condition (\ref{eq:A(brzeg) zawarte w R^u minus B(0,1)}) implies that $%
\left\| A_{i}p\right\| >\left\| p\right\| $ for $i=1,\ldots n+1$ and $p\neq
0 $, which gives us $p_{0}=0$. The fact that $p_{1}=\ldots =p_{n}=0$ follows
from (\ref{eq:pi=Api-1}). \qed

\bigskip

\noindent \textbf{Proof of Lemma \ref{lem:lemat dla lambda in [0,1]}:} From
the homotopy property, it is sufficient to show that
\begin{equation}
H(\lambda ,x)\neq 0,\quad \text{for all }x\in \partial XN\text{ and }\lambda
\in \lbrack 0,1].  \label{eq:H rozne od zera na
brzegu}
\end{equation}

We will consider an $x$ from the boundary of $XN$ \newline
$x=((x_{0},t_{0}),\ldots ,(x_{n},t_{n}))$. If $x\in \partial XN$ then there
exists an $i$ such that one of the following conditions holds
\begin{eqnarray}
x_{i} &\in &N_{i}^{+}  \label{eq:xi in Ni+} \\
x_{i} &\in &N_{i}^{-}  \label{eq:xi in Ni-} \\
t_{i} &\in &\{s_{i}-r_{i},s_{i}+r_{i}\}  \label{eq:si+ri}
\end{eqnarray}

First let us consider the case (\ref{eq:xi in Ni+}). For $i=1,\ldots n$ if $%
x_{i}\in N_{i}^{+}$ and $H(\lambda ,x)=0$ then in particular
\begin{equation}
x_{i}-H^{i}(\lambda ,t_{i-1},x_{i-1})=0  \label{eq:warunek dla N+}
\end{equation}
From the statement of Lemma \ref{lem:Lemma o nakrywaniu dla perturbacji},
condition (\ref{eq:H4}) we know that
\begin{equation*}
H^{i}([0,1],t_{i-1},N_{i-1})\cap N_{i}^{+}=\emptyset .
\end{equation*}
This and the fact that $x_{i-1}\in N_{i-1}$ contradicts
(\ref{eq:warunek dla N+}). We therefore know that (\ref{eq:xi in
Ni+}) does not hold for $i=1,\ldots,n$. For $i=0$ if $x_{0}\in
N_{0}^{+}$ and $H(\lambda ,x)=0$ then
\begin{eqnarray*}
x_{0}-G(\lambda ,t_{n},x_{n}) &=&0 \\
(p_{0},q_{0})-(A_{n+1},0) &=&0
\end{eqnarray*}
which means that $q_{0}=0$ which contradicts the fact that $%
x_{0}=(p_{0},q_{0})\in N_{0}^{+}=\overline{B_{u}}(0,1)\times \partial
\overline{B_{s}}(0,1).$

Let us now consider the case (\ref{eq:xi in Ni-}). For $i=0,\ldots ,n-1$ if $%
x_{i}\in N_{i}^{-}$ and $H(\lambda ,x)=0$ then
\begin{equation}
x_{i+1}-H^{i+1}(\lambda ,t_{i},x_{i})=0  \label{eq:warunek dla N-}
\end{equation}
From Lemma \ref{lem:Lemma o nakrywaniu dla perturbacji}, condition (\ref
{eq:H3}) we have
\begin{equation*}
H^{i+1}([0,1],t_{i},N_{i}^{-})\cap N_{i+1}=\emptyset ,
\end{equation*}
which contradicts (\ref{eq:warunek dla N-}). Condition (\ref{eq:xi in Ni-}) cannot hold for $i=0,\ldots,n-1$. For $i=n$
if $x_{n}\in N_{n}^{-}$ and $%
H(\lambda ,x)=0$ then

\begin{eqnarray*}
x_{0}-G(\lambda ,t_{n},x_{n}) &=&0 \\
(p_{0},q_{0})-(A_{n+1}p_{n},0) &=&0
\end{eqnarray*}
The fact that $x_{n}\in N_{n}^{-}$ means that $p_{n}\in \partial \overline{%
B_{u}}(0,1).$ We know that $p_{0}\in $ $\overline{B_{u}}(0,1)$ and $%
p_{0}=A_{n+1}p_{n},$ which contradicts the fact that $A_{n+1}(\partial
\overline{B_{u}}(0,1))\subset R^{u}\backslash \overline{B_{u}}(0,1).$

We are now left with the case (\ref{eq:si+ri}). For $i=1,\ldots ,n$ if (\ref
{eq:si+ri}) holds and $H(\lambda ,x)=0$ then in particular
\begin{equation}
t_{i}-\lambda s_{i}-(1-\lambda )\sigma_i (x_{i-1},t_{i-1})=0.
\label{eq:warunek dla ti}
\end{equation}
Our construction of $XN$ (\ref{eq:konstrukcja XN}) which guarantees that
\begin{equation*}
F_{i}(N_{i-1},I_{i-1})\subset V_{i}\times \text{int}I_{i}.
\end{equation*}
gives us
\begin{eqnarray*}
\sigma_i (x_{i-1},t_{i-1}) &\in &(s_{i}-r_{i},s_{i}+r_{i}) \\
\lambda s_{i}+(1-\lambda )\sigma_i (x_{i-1},t_{i-1}) &\in
&(s_{i}-r_{i},s_{i}+r_{i})
\end{eqnarray*}
and therefore from (\ref{eq:si+ri})
\begin{equation*}
t_{i}-\lambda s_{i}-(1-\lambda )\sigma_i (x_{i-1},t_{i-1})\neq 0
\end{equation*}
This clearly contradicts (\ref{eq:warunek dla ti}). For $i=0$ from the
definition of $H(\lambda ,\cdot )$ and the fact that $H(\lambda ,x)=0$ we
get straight away the fact that $t_{0}=0$ which means that it is not
possible for $t_{0}\in \{-1,1\}.$

We have shown that for any $\lambda \in \lbrack 0,1]$ and $x\in \partial XN$%
, $H(\lambda ,x)\neq 0,$ this fact and the homotopy property of the index
concludes our proof. \qed

\bigskip

\textbf{Proof of Theorem \ref{thm:basicper}:} Let us consider the sequence
of covering relations from $\Gamma $
\begin{equation}
N_{0}\cover{f_{01}}N_{1}\cover{f_{12}}N_{2}\cover{f_{23}}\ldots
\label{eq: ciag relacji z Gamma}
\end{equation}
Let us consider a finite subsequence of the sequence (\ref{eq:
ciag relacji z Gamma})

\begin{equation*}
N_{0}\cover{f_{01}}N_{1}\cover{f_{12}}\ldots \cover{f_{m-1\, m}}N_{m}
\end{equation*}
From Lemma \ref{lem: lemat glowny} we know that for $\left| \epsilon \right|
<\delta $ there exists \bigskip $x_{m}\in N_{0}$ and a sequence $%
t_0=t_{0}^m<t_{1}^m<\ldots <t_{m}^m$ , such that for the solution
$x$ of the equation (\ref{eq:x'=v+e}) we have
\begin{equation*}
x(t_{i})\in \text{int}N_{i}\quad \text{for }i=1,\ldots ,m
\end{equation*}
and that $\delta $ depends only on the family $\Gamma $ and not on the
length of the sequence. We therefore have a sequence $\{x_{m}\}_{m=1}^{%
\infty }\subset N_{0}$. Since $N_{0}$ is compact there exists a subsequence $%
x_{m_{k}}$ which converges to a certain $x_{0}\in N_{0}$. In the curse of
the proof of Lemma \ref{lem: lemat glowny} we have shown that (\ref
{eq:ti-sigma=0})
\begin{equation*}
t_{i}^{m}-\sigma_i (x(t_{i-1}^{m}),t_{i-1}^{m})=0
\end{equation*}
which together with the fact from (\ref{eq:sigma-T ograniczone}), that $%
\sigma_i (x(t_{i-1}^{m}),t_{i-1}^{m})-t_{i-1}^m$ is bounded, means that $%
t_{i}^m-t_{i-1}^m=\sigma_i (x(t_{i-1}^{m}),t_{i-1}^{m})-t_{i-1}^m$ is
bounded. From this fact and from the continuity of the solution of the
problem
\begin{equation*}
x^{\prime }(t)=v(x(t))+\epsilon (t,x)
\end{equation*}
with respect to the initial conditions, it follows that the solution $x(t),$
of the problem
\begin{eqnarray}
x^{\prime } &=&v(x)+\epsilon (t,x) \\
x(0) &=&x_{0}  \notag
\end{eqnarray}
passes through the sets $N_{0},N_{1},\ldots $ and therefore there exists a
sequence $t_{0}<t_{1}<\ldots $ such that
\begin{equation*}
x(t_{i})\in \text{int}N_{i}\quad \text{for }i=1,2,\ldots
\end{equation*}
\qed

\section{Application to R\"ossler equations.}
\label{sec:appl} In this section we combine
Theorem~\ref{thm:basicper} and  results from \cite{ZNon} to show
that small non-autonomous perturbations of R\"{o}ssler \cite{R} 
posses symbolic dynamics.

First we need to recall some definitions.

Let $k$ be a positive integer. Let $\Sigma_k
:=\{0,1,\dots,k-1\}^\mathbb{Z}$, $\Sigma^+_k
:=\{0,1,\dots,k-1\}^\mathbb{N}$.  $\Sigma_k$, $\Sigma^{+}_k$ are
topological spaces with the Tichonov topology. On $\Sigma_k$,
$\Sigma^{+}_k$ we have the shift map $\sigma$ given by
\begin{displaymath}
  (\sigma(c))_i=c_{i+1}
\end{displaymath}

Let $A = [\alpha_{ij}] $ be a $k \times k$-matrix, \newline
$\alpha_{ij} \in \mathbb{R}_+ \cup \{0\}$, $i,j = 0,1,\dots,k-1$.
We define   $\Sigma_A \subset \Sigma_k$ and $\Sigma^+_A \subset
\Sigma^+_k$ by
\begin{eqnarray}
   \Sigma_A & := \{ c=(c_i)_{i \in \mathbb{Z}} \quad | \quad
       \alpha_{c_i c_{i+1}} > 0 \} \\
   \Sigma^{+}_A & := \{ c=(c_i)_{i \in \mathbb{N}} \quad | \quad
       \alpha_{c_i c_{i+1}} > 0 \}
\end{eqnarray}
Obviously $\Sigma^+_A$, $\Sigma_A$ are invariant under $\sigma$.

Let $F:X \to X$ be any continuous map  and  $N \subset X$. By
$F_{|N}$ we will denote the map obtained by  restricting  the
domain of $F$ to the set $N$. The maximal invariant part of $N$
(with respect to $F$) is defined by
\begin{displaymath}
       \mbox{Inv}(N,F) = \bigcap_{i \in \mathbb{Z}} F_{|N}^{-i}(N).
\end{displaymath}

The R\"ossler equations are given by \cite{R}
\begin{eqnarray}
   \dot x &  = & -(y+z) \nonumber \\
   \dot y &  = & x + by  \label{eq:rossler}  \\
   \dot z &  = & b + z(x-a)  \nonumber
\end{eqnarray}
where $a=5.7, b=0.2$. These are parameters values originally
considered by R\"ossler. The  flow generated  by Eq.
(\ref{eq:rossler}) exhibits a so-called strange attractor.

We will investigate  the Poincar\'e map $P$  generated  by
(\ref{eq:rossler}) on the section $\Theta:=\{(x,y,z)| \quad x=0, y
<0, \dot x > 0 \}$.

The following result was proved in \cite{ZNon} (see also
\cite{Zbif})
\begin{theorem}
\label{thm:ross}

For all parameter values in  sufficiently small neighborhood of
$(a,b)=(5.7,0.2)$ there exists Poincar\'e section $N \subset
\Theta$ such that the Poincar\'e map $P$ induced by Eq.
(~\ref{eq:rossler}) is well defined and continuous.

There exists continuous map $ \pi : \mbox{Inv}(N,P) \to \Sigma_3
$, such that
\begin{displaymath}
  \pi \circ P = \sigma \circ \pi.
\end{displaymath}
$\Sigma_A \subset \pi(\mbox{Inv}(N,P))$, where
\begin{displaymath}
  A:= \left[
    \begin{array}{ccc}
     0 & 1 & 1 \\
     0 & 1 & 1 \\
     1 & 0 & 0
    \end{array}
    \right]
\end{displaymath}

The preimage of any periodic sequence from $\Sigma_A$ contains
periodic points of $P$.
\end{theorem}

Above theorem is a consequence of Theorem~\ref{thm:t1} and the
following Lemma, which was established in \cite{ZNon} with
computer assistance (computer assisted proof)
\begin{lemma}
\label{lem:pom} There are $h$-sets $N_0,N_1,N_2 \subset \Theta$
such that for all parameter values in sufficiently small
neighborhood of $(a,b)=(5.7,0.2)$ $N \subset \mbox{Dom}(P)$ and
the following conditions hold
\begin{equation}
  N_0 \cover{P} N_2, \quad
  N_1 \cover{P} N_0, N_1, \quad
  N_2 \cover{P} N_0, N_1
\end{equation}
\end{lemma}

Let us denote by $R_{a,b}:\mathbb{R}^3 \to \mathbb{R}^3$ the
vector field on the right-hand side of (\ref{eq:rossler}). By
applying Theorem~\ref{thm:basicper} to Lemma~\ref{lem:pom} we
immediately obtain the following
\begin{theorem}
\label{thm:rossper} Let us fix $(a,b)=(5.7, 0.2)$. Let $A$ be as
in Theorem~\ref{thm:ross}. Consider a non-autonomous perturbation
of (\ref{eq:rossler})
\begin{equation}
  x'=R_{a,b}(x) + \epsilon(t,x). \label{eq:rossper}
\end{equation}

 There exists $\delta>0$, such that for any $t_0 \in \mathbb{R}$
 and any sequence $c=(c_i) \in \Sigma_A^+$ there exists a solution
 of (\ref{eq:rossper}), $x_c:[t_0,\infty) \to \mathbb{R}^3$ and a
 sequence $t_0 < t_1 <  t_2 < \dots < t_n < t_{n+1} < \dots $,
 such that
 \begin{eqnarray*}
   x_c(t) &\in& \Theta, \qquad \mbox{iff  $t=t_i$ for some i } \\
   x_c(t_i) &\in& |N_{c_i}|.
 \end{eqnarray*}
\end{theorem}

Above theorem says nothing about the size of $\delta$. To obtain a
numerical value for  $\delta$ one can take one of two approaches
\begin{description}
\item[analytical] from the computer assisted proof in \cite{ZNon}
one can obtain global bounds $Z \subset \mathbb{R}^3$, $Z$
compact, such that all trajectories linking $|N_i|$ with its
Poincar\'e image are in $Z$. For $\epsilon$ sufficiently small the
same will be true for (\ref{eq:rossper}). Now using bounds for the
Poincar\'e return times on $|N_i|$ we can compute an upper bound
of the distance between the solution of (\ref{eq:rossler}) and
(\ref{eq:rossper}). Then we compute $\epsilon$ for which the
covering relations listed in Lemma~\ref{lem:pom} survive.
\item[computational] we can replace (\ref{eq:rossper}) by a
differential inclusion
\begin{equation}
  x' \in R_{a,b}(x) + [-\delta,\delta]^3  \label{eq:rossdiffincl}.
\end{equation}
Now for various values of $\delta$ we can perform an rigorous
integration of (\ref{eq:rossdiffincl}) looking for the largest
possible $\delta$   for which the covering relations listed
Lemma~\ref{lem:pom} are still satisfied (for any continuous
selector). For an algorithm for rigorous integration of
differential inclusions see  \cite{ZKSPer}.
\end{description}

\subsection{Other examples.}
Other  based on Theorem~\ref{thm:t1} computer assisted proofs of
the existence of nontrival symbolic dynamics give rise to theorems
analogous to Theorem~\ref{thm:rossper}. This applies to the
following systems
\begin{itemize}
\item Lorenz equations  \cite{GaZ},
\item Chua circuit \cite{G},
\item Kuramoto-Shivashinsky ODE \cite{W}.
\end{itemize}
The precise formulation of these results is left to the reader.

\section{ Appendix. Properties of the local Brouwer degree}
\label{sec:app-br-deg}

\bigskip \textbf{Homotopy property.} \cite{L} Let $H:[0,1]\times
D\rightarrow R^{n}$ be continuous. Suppose that
\begin{equation}
\bigcup_{\lambda \in \lbrack 0,1]}H_{\lambda }^{-1}(c)\cap D\quad
\text{is compact}  \label{eq:union}
\end{equation}
then
\begin{equation*}
\forall \lambda \in \lbrack 0,1]\quad \text{deg}(H_{\lambda },D,c)=\text{deg}%
(H_{0},D,c)
\end{equation*}
If $[0,1]\times \overline{D}\subset $dom$(H)$ and $\overline{D}$
is compact, then (\ref{eq:union}) follows from the condition
\begin{equation*}
c\notin H([0,1],\partial D)
\end{equation*}

\textbf{Degree property for affine maps.} \cite{L} Suppose that $%
f(x)=B(x-x_{0})+c$, where $B$ is a linear map and $x_{0}\in
R^{n}.$ If the equation $B(x)=0$ has no nontrivial solutions (i.e
if $Bx=0$, then $x=0$) and $x_{0}\in D$, then
\begin{equation}
\text{deg}(f,D,C)=\text{sgn}(\text{det}B).  \label{eq:deg(f,C)
=sgn(det(B))}
\end{equation}

\end{document}